\documentclass[11pt,reqno]{article}
\usepackage{latexsym}
\usepackage{amssymb,amsmath}
\usepackage{graphics}
\usepackage{epsfig}
\usepackage[]{times}
\makeatletter \def\section{\@startsection {section}{1}{\z@}{-3.5ex plus
-1ex minus -.2ex}{2.3ex plus .2ex}{\large\sc}}
\def\subsection{\@startsection{subsection}{2}{\z@}{-3.25ex plus -1ex minus
 -.2ex}{1.5ex plus .2ex}{\normalsize\sc}}
\makeatother
\makeatletter
\@addtoreset{equation}{section}

\makeatother
\newcommand{\nc}{\newcommand}

\nc{\bea}{\begin{eqnarray}} \nc{\eea}{\end{eqnarray}}
\nc{\be}{\bea} \nc{\ee}{\eea}
\nc{\tr}{\mathop{\mbox{tr}}\nolimits}
\nc{\ad}{\mathop{\mbox{ad}}\nolimits}
\nc{\Tr}{\mathop{\mbox{Tr}}\nolimits}
\nc{\Det}{\mathop{\mbox{Det}}\nolimits}
\nc{\rk}{\mathop{\mbox{rk}}\nolimits}
\nc{\ra}{\rightarrow}
\nc{\Ra}{\Rightarrow}
\nc{\LRa}{\Leftrightarrow}
\nc{\ot}{\otimes}
\nc{\non}{\nonumber\\}
\nc{\ZZ}{\mathbb{Z}}
\nc{\RR}{\mathbb{R}}
\newtheorem{theorem}{Theorem}

\newenvironment{proof}[1][Proof]{\begin{trivlist}
\item[\hskip\labelsep {\bfseries #1}]}{\end{trivlist}}
\def\p{\prime}
\def\1#1{^{(#1)}}
\def\la{\langle}
\def\ra{\rangle}
\def\be{\begin{equation}}
\def\ee{\end{equation}}
\def\bea{\begin{eqnarray}}
\def\eea{\end{eqnarray}}
\begin{document}
\title{Wilf's question in numerical semigroups $S_3$ revisited
\\and inequalities for rescaled genera}
\author{Leonid G. Fel\\ \\
Department of Civil Engineering, Technion - Israel Institute of
Technology,\\
Haifa 32000, Israel\\
{\em e-mail: lfel@technion.ac.il}}
\date{}
\maketitle
\vspace{-1cm}
\def\be{\begin{equation}}
\def\bea{\begin{eqnarray}}
\def\ee{\end{equation}}
\def\eea{\end{eqnarray}}
\def\p{\prime}
%%%%%%%%%%%%%%%%%%%%%%%%%%%%%%%%%%%%%%%%
\begin{abstract}
We consider numerical semigroups $S_3=\la d_1,d_2,d_3
\ra$, minimally generated by three positive integers. We revisit the Wilf question in $S_3$ and, making use of identities for degrees of syzygies of such 
semigroups, give a short proof of existence of an affirmative  answer. We find the  upper and lower bounds for rescaled genera of numerical semigroups 
$S_3$.\\\\
\noindent
{\bf Key words}: numerical semigroups, identities for degrees of 
syzygies, Wilf's question, \\
{\bf 2000 Math. Subject Classification:} Primary -- 20M14, Secondary -- 11P81.
\end{abstract}
%%%%%%%%%%%%%%%%%%%%%%%%%%%%%%%%%%%%%%%%%%%%%%%%%%%%%%%%%%%%%%%%
%%%%%%%%%%%%%%%%%%%%%%%%%%%%%%%%%%%
\section{Introduction}\label{l1}
%%%%%%%%%%%%%%%%%%%%%%%%%%%%%%%%%%%%%%%%%%%%%%%%%
Let a numerical semigroup $S_3$ be minimally generated by a set of natural numbers 
$\{d_1,d_2,d_3\}$, where
\bea
\gcd(d_1,d_2,d_3)=1,\qquad 3\le d_1<d_2<d_3\le d_1d_2-d_1-d_2,\label{a1}
\eea
and neither of its generators is linearly representable by the rest of them. Its generating function $H\left(S_3;t\right)$
\bea
H\left(S_3;t\right)=\sum_{s\;\in\;S_3}t^s,\qquad t<1,\qquad 0\in S_3,\label{a2}
\eea
is referred to as {\em the Hilbert series} of $S_3$ and has a rational representation (Rep),
\bea
H\left(S_3;t\right)=\frac{1-t^{x_1}-t^{x_2}-t^{x_3}+t^{y_1}+t^{y_2}}{\left(1-t^{d_1}\right)
\left(1-t^{d_2}\right)\left(1-t^{d_3}\right)},\qquad x_j,y_j\in {\mathbb Z}_>,\label{a3}
\eea
where $x_j,y_j$ denote degrees of {\em syzygies}. The Frobenius number $F_3$ as a largest unrepresentable integer by the triple $\{d_1,d_2,d_3\}$, is often used in commutative algebra when shifting by unity, $c_3=F_3+1$, and referred to as {\em a conductor} of $S_3$. The largest degree $g_3=\max\left\{y_1,y_2\right\}$ is related to $F_3$,
\bea
g_3=F_3+\sigma_1,\quad \sigma_1=d_1+d_2+d_3,\qquad G_0=\#\Delta_3,\quad
\Delta_3={\mathbb Z}_>\setminus S_3,\label{a4}
\eea
where $\Delta_3$ and $G_0$ stand for the set of unrepresentable inetgers ({\em gaps}) and {\em the genus} of semigroup $S_3$. 

The set $\Delta_3$ comprises two sorts of gaps: those $s\in \Delta_3$ such that $F_3-s\not\in\Delta_3$, and those $s\in \Delta_3$ such that $F_3-s\in\Delta_3$. Such dichotomy of the set $\Delta_3$ leads to a simple inequality
\bea
w_3\ge\frac1{2},\qquad w_3=\frac{G_0}{c_3},\label{a5}
\eea
and equality in (\ref{a5}) is arisen only for symmetric semigroups where $\Delta_3$ comprises only the 1st sort gaps. In fact, inequality (\ref{a5}) 
holds for numerical semigroups $S_m$ of any $m$, i.e., $w_m\!\ge\! 1/2$ where $w_m\!=\!G_0/c_m$.

In 1978, Wilf \cite{wl78} raised two questions ({\sf WQ}) and the 1st of them was: is it true that for given $S_m$ holds
\bea
w_m\le\frac{m-1}{m},\label{a6}
\eea
with equality only for the generators $m,m+1,\ldots,2m-1$ ? In a seminal paper \cite{fr87}, by means of embedding procedure for a sequence of partial gapsets of semigroup $S_m$, it was shown
\bea
w_m\le\frac{\tau_m}{\tau_m+1},\label{a7}
\eea
where $\tau_m$ stands for the {\em type} of $S_m$. Inequality (\ref{a7}) coincides with (\ref{a6}) if $m=2,3$, but does not imply (\ref{a6}) if $m>3$. During the last decade, a vast literature \cite{de20} was devoted to Wilf's question which, despite the attention it attracted in various special cases, remains unsolved completely for non-symmetric semigroups $S_m$, $m\ge 4$.

Based on polynomial identities for degrees of syzygies \cite{fe17}, in the present paper, we suggest quite different approach to get an 
affirmative answer to {\sf WQ} in semigroups $S_3$. In  its framework we obtain also the lower bound for $F_3$ as well as the lower and upper bounds for rescaled genera of $S_3$. From the standpoint of {\sf 
WQ} perspectives, the developed approach may be applied to semigroups $S_m$, $m\ge 4$. We plan to address this question in a separate paper.
%%%%%%%%%%%%%%%%%%%%%%%%%%%%%%%%%%%%%%%%%%%%%%%%%
\section{Syzygy identities for numerical semigroups $S_3$}\label{l2}
%%%%%%%%%%%%%%%%%%%%%%%%%%%%%%%%%%%%%%%%%%%%%%%%%
Consider a numerical semigroup $S_3$ and write the polynomial identities (see \cite{fe17}, Theorem 1) for degrees of syzygies $x_j,y_j$. Denote by $X_r$ and $Y_r$ two power sums,
\bea
a)\quad X_r=x_1^r+x_2^r+x_3^r,\qquad b)\quad Y_r=y^r+g_3^r,\qquad 0\le x_j,y<g_3,\label{b1}
\eea
and consider a set of polynomial equations for five real variables $x_1,x_2,x_3,y$ and $g_3$,
\bea
&&a)\quad Y_1-X_1=0,\nonumber\\
&&b)\quad Y_2-X_2=2\pi_3,\qquad
\pi_3=d_1d_2d_3,\nonumber\\
&&c)\quad Y_{r+3}-X_{r+3}=\frac{(r+3)!}{r!}\;
K_r\pi_3,\quad r\ge 0,\quad K_r>0,\label{b2}
\eea
where coefficients $K_r$ exhibit a linear combination of higher genera $G_0,\ldots,G_r$ of numerical semigroup $S_3$, i.e., $G_r=\sum_{s\in\Delta_3}\!s^r
$, $r\ge 0$ (see formulas (22,23) in \cite{fe21}). 

E.g.,
\bea
\hspace{-.5cm}
&&K_0=G_0+\delta_1,\hspace{6.5cm}\delta_k=\frac{\sigma_k-1}{2^k},\label{b3}\\
&&K_1=G_1+\frac{\sigma_1}{2}G_0+\frac{3\delta_1^2+\delta_2}{6},\hspace{4cm}
\sigma_k=d_1^k+d_2^k+d_3^k,\nonumber\\
&&K_2=G_2+\sigma_1G_1+\frac{3\sigma_1^2+\sigma_2}{12}G_0+\frac{\delta_1(\delta_1^2+\delta_2)}{3},
\hspace{1.2cm}\mbox{etc.}\nonumber
\eea
Explicit expressions of $K_r$ are given in \cite{fe21}, formula (27).
%%%%%%%%%%%%%%%%%%%%%%%%%%%%%%%%%%%%%%%%%%%%%%%%%
\subsection{Lower bound of $g_3$}\label{l21}
%%%%%%%%%%%%%%%%%%%%%%%%%%%%%%%%%%%%%%%%%%%%%%%%%
Make use of Newton-Maclaurin 's inequalities \cite{har59} for power sums $X_j$ and $Y_j$,
\bea
a)\quad 3X_2>X_1^2,\qquad b)\quad 2Y_2>Y_1^2,\label{b4}
\eea
and substitute equality (\ref{b2}b) into inequality (\ref{b4}a)
\bea
3\left(Y_2-2\;\pi_3\right)>Y_1^2.\label{b5}
\eea
By (\ref{b1}b), find a relationship between $Y_1$ and $Y_2$,
\bea
Y_2=Y_1^2-2g_3Y_1+2g_3^2,\label{b6}
\eea
and substitute it into (\ref{b5}),
\bea
Y_1^2-3g_3Y_1+3g_3^2>3\pi_3,\qquad g_3<Y_1<2g_3.\label{b7}
\eea
Denote by $u$ a ratio $Y_1/g_3$ and rescale (\ref{b7}) by $g_3$ as follows,
\bea
\frac{g_3^2}{3\pi_3}>P(u),\qquad P(u)=\frac1{u^2-3u+3},\qquad 
u=\frac{Y_1}{g_3},\qquad 1<u<2.\label{b8}
\eea
Find a range of $g_3$ where inequality (\ref{b8}) is satisfied for any $u\in(1,2)$. Since a convex function $P(u)$ reached its minimum $P(u)=1$ at $u=1,2$, we arrive at Davison's lower bound \cite{da94} for $g_3$,
\bea
g_3>\sqrt{3\;\pi_3}.\label{b9}
\eea
In fact, more accurate reasoning \cite{fe06}, e.g., $Y_1\le 2g_3-1$, leads to a slightly stronger bound, $\sqrt{3}\sqrt{\pi_3+1}$.
%%%%%%%%%%%%%%%%%%%%%%%%%%%%%%%%%%%%%%%%%%%%%%%%%
\subsection{Lower and upper bounds of $h_0=K_0/g_3$}\label{l22}
%%%%%%%%%%%%%%%%%%%%%%%%%%%%%%%%%%%%
Introduce elementary symmetric polynomials
\bea
{\cal X}_1=x_1+x_2+x_3,\quad {\cal X}_2=
x_1x_2+x_2x_3+x_3x_1,\quad
{\cal X}_3=x_1x_2x_3,\label{b10}
\eea
which are related by Newton's recursion identities to power sums $X_k$, defined in (\ref{b1}a),
\bea
{\cal X}_1=X_1,\qquad 2{\cal X}_2=X_1^2-X_2,\qquad
6{\cal X}_3=X_1^3-3X_1X_2+2X_3.\label{b11}
\eea
Recall the Newton-Maclaurin inequalities \cite{har59} for polynomials ${\cal X}_k$,
\bea
a)\quad\frac{{\cal X}_1}{3}\ge\left(\frac{{\cal X}_2}{3}\right)^{1/2}\ge {\cal X}_3^{1/3},\qquad 
b)\quad\left(\frac{{\cal X}_2}{3}\right)^2\ge
\frac{{\cal X}_1}{3}{\cal X}_3.\label{b12}
\eea
Consider the 2nd inequality in (\ref{b12}a), and substitute there identiies (\ref{b11})
\bea
(X_1^2-X_2)^3\ge 6(X_1^3-3X_1X_2+2X_3)^2.\nonumber
\eea
Now, substitute into the last inequality three first equalities (\ref{b2} a,b,c),
\bea
(Y_1^2-Y_2+2\pi_3)^3\ge 6(Y_1^3-3Y_1Y_2+
6Y_1\pi_3+2Y_3-12\pi_3K_0)^2.\label{b13}
\eea
Combining identity for power sums, $2Y_3=3Y_1Y_2-Y_1^3$, and relation (\ref{b6}), we simplify (\ref{b13})
\bea
(g_3Y_1-g_3^2+\pi_3)^3\ge 27\pi_3^2\;(2K_0-Y_1)^2.\label{b14}
\eea
Introduce two new variables $v$ and $h_0$,
\bea
v=\frac{\pi_3}{g_3^2},\quad 0<v<\frac1{3},\qquad h_0=\frac{K_0}{g_3},\quad\frac1{2}\le h_0\le w_3,
\label{b15}
\eea
and rescale (\ref{b14}) by $g_3$
\bea
(u-2h_0)^2\le Q(u,v),\qquad
Q(u,v)=\frac{(u-1+v)^3}{27v^2}<Q(2,v).\label{b16}
\eea
In (\ref{b15}), $v$ is bounded from above due to (\ref{b9}) and $h_0$ is bounded from below due to (\ref{a5}).

Find the upper bound of $h_0$ and represent (\ref{b16}) as follows,
\bea
u<2h_0+\sqrt{Q(2,v)}.\label{b17}
\eea
On the other hand, there holds always another inequality, $u<2$. In order to find the maximal value of $h_0$ when both inequalities still hold, we have to choose such $v$, which provides minimal value of $Q(2,v)$ at the interval $v\in(0,1/3)$, that happens at $v=1/3$.
%%%%%%%%%%%%%%%%%%%%%%%%%%%%%%%%
\begin{figure}[h!]\begin{center}
\psfig{figure=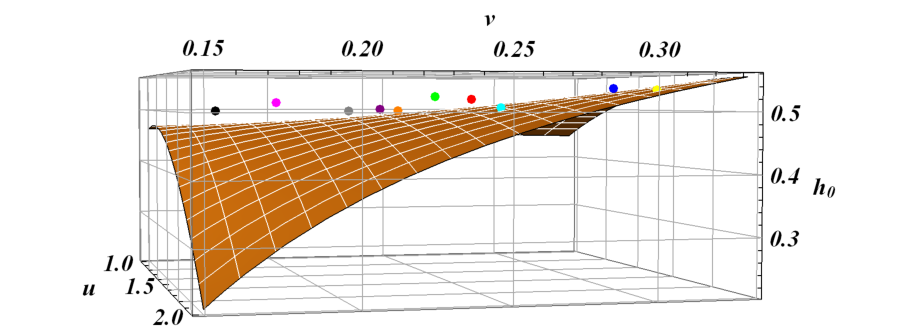,height=6cm}\\
\end{center}
\vspace{-.5cm}
\caption{Plot of the function $
\Phi(u,v)\!\!=\!\!\frac1{2}\left(u\!-\sqrt{Q(u,v)}\right)$ and the points $(u,v,h_0)$ for semigroups 
$\la 3,4,5\ra$  ({\it yellow}), $\la 4,5,6\ra$  
({\it purple}), $\la 5,6,7\ra$ ({\it blue}), $\la 11,17,29\ra$ ({\it red}), $\la 25,31,43\ra$ ({\it magenta}), $\la 23,29,44\ra$ ({\it cyan}), 
$\la 43,47,113\ra$ ({\it green}), $\la 501,503,603\ra$ ({\it black}), $\la 901,903,1003\ra$ ({\it gray}) and $\la 1201,1203,1303\ra$ ({\it orange}).}\label{fig1}
\end{figure}
%%%%%%%%%%%%%%%%%%%%%%%%%%%%%%%%%%%%
Thus, we arrive at equality
\bea
2h_0<2-\sqrt{Q(2,1/3)}=\frac{10}{9},\label{b18}
\eea
and finally (see Figure \ref{fig1}),
\bea
\frac1{2}\le h_0<\frac{5}{9}\simeq 0.5555.\label{b19}
\eea
%%%%%%%%%%%%%%%%%%%%%%%%%%%%%%%%%%%%%%%%%%%%%%%%%
\subsection{Wilf's question in $S_3$.}\label{l23}
%%%%%%%%%%%%%%%%%%%%%%%%%%%%%%%%%%%%%%%%%%%%%%%%%
Prove that {\sf WQ} has an affirmative answer for all non-symmetric semigroups $S_3$ and start with $\la 3,d_2,d_3\ra$. Making use of Lemma 6 in [8], which established inequalities for such semigroups, we get: 
$$
\frac{3c_3}{2}\le 3G_0<2c_3+1.
$$ 
Replacing the r.h.s. of strict inequality on non-strict one, $3G_0\le 2c_3$, we arrive at $w_3\le 2/3$.

Consider numerical semigroups $S_3$ with $d_1\ge 4$ and, using relation $g_3=c_3+2\delta_1$, and definition (\ref{b15}) of $h_0$ and expression 
(\ref{b3}) for $K_0$, represent inequality (\ref{b19}) for the upper bound as follows,
\bea
w_3<\frac{5}{9}+\frac{e_3}{9},\qquad e_3=\frac{\delta_1}{c_3}.\label{b20ā}
\eea
A sufficient (not necessary) condition for affirmative answer to {\sf WQ} would be inequality $e_3\le 1$. 
%%%%%%%%%%%%%%%%%%%%%%%%%%%%%%%%%%%%%%%%%%%%%%%%%%%%%%%%%%%%
\begin{theorem}\label{the1}
Let a non-symmetric numerical semigroup $\la d_1,
d_2,d_3\ra$, $d_1\ge 4$, be given and its generators satisfy (\ref{a1}). The {\sf WQ} has an affirmative answer for all numerical semigroups $S_3$ 
\end{theorem}
%%%%%%%%%%%%%%%%%%%%%%%%%%%%%%%%%%%%%%%%%%%%%%%%%%%%%%%%%%%%
\begin{proof}
Instead of $e_3$, consider its inverse $1/e_3$ and apply inequality (\ref{b9}),
\bea
\frac1{e_3}=\frac{g_3}{\delta_1}-2>
2\left(\frac{\sqrt{3\pi_3}}{\sigma_1-1}-1\right)>
2\left(\frac{\sqrt{3\pi_3}}{\sigma_1}-1\right).\label{b21}
\eea
According to (\ref{b21}), a sufficient (not necessary) condition to provide $e_3<1$ is the Diophantine inequality
\bea
\rho_3<\frac{2}{\sqrt{3}}\simeq 1.1547,\qquad 
d_1\ge 4,\qquad \rho_3=\frac{d_1+d_2+d_3}
{\sqrt{d_1d_2d_3}}.\label{b22}
\eea
Check its solvability for different triples $\{d_1,
d_2,d_3\}$. For that purpose, represent (\ref{b22}) as follows
\bea
(\sqrt{d_3}-\sqrt{d_2})^2<2\left(\sqrt{\frac{d_1}{3}}-1\right)\sqrt{d_2d_3}-d_1,\label{b23}
\eea
that necessary leads to the lower bound of the product $d_2d_3$, 
\bea
d_2d_3>C(d_1),\qquad C(d_1)=\frac{3}{4}
\left(\frac{d_1}{\sqrt{d_1}-\sqrt{3}}\right)^2.\label{b24}
\eea
where the concave function $C(d_1)$ arrives its minimum $C(12)=36$ (see Figure \ref{fig2}).
%%%%%%%%%%%%%%%%%%%%%%%%%%%%%%%%%%%%%%%%%%%
\begin{figure}[h!]
\begin{center}
\psfig{figure=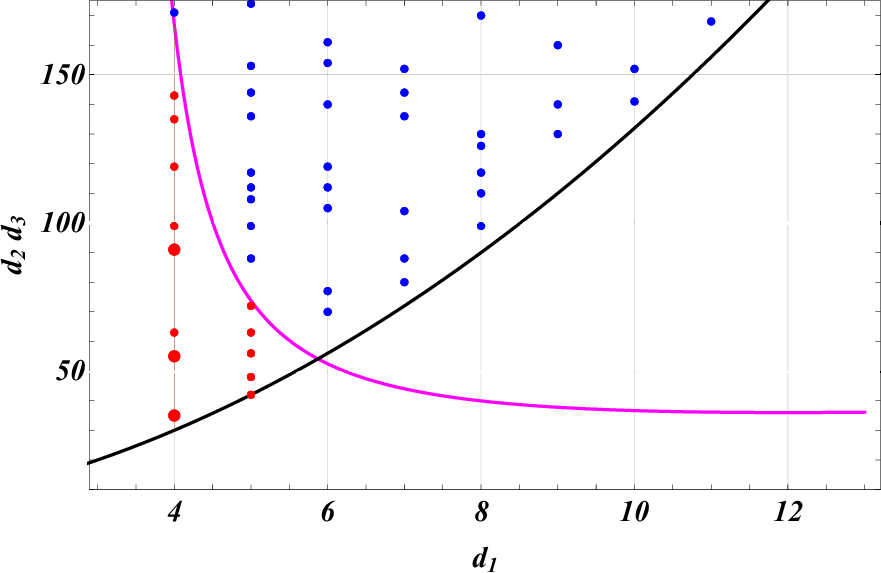,
height=6.5cm}\\
\end{center}
\vspace{-.1cm}
\caption{Plot of functions $C(d_1)$ ({\em in magenta}) and $(d_1+1)(d_1+2)$ ({\em in black}). 
{\em Blue dots} correspond to semigroups with 
$\rho_3<1.1547$, $e_3<1$. Ten small and three large 
{\em red dots} correspond to semigroups with 
$\rho_3>1.1547$, $e_3<1$ and $\rho_3>1.1547$, $e_3>1$, respectively. These semigroups are listed in 
Table 1.}\label{fig2}
\end{figure}
%%%%%%%%%%%%%%%%%%

The criterion (\ref{b24}) has to be supplemented by another restriction, $d_3>d_2>d_1$, or
\bea
d_2d_3\ge(d_1+1)(d_1+2),\label{e25}
\eea

The first values of descendent sequence C(d1) are given below,
\bea
C(4)=167.14,\quad C(5)=73.81,\quad C(6)=52.46,\quad C(7)=44.02,\quad\ldots.\nonumber
\eea
The final criterion for generators $d_j$ to satisfy (\ref{b24}) reads,
\bea
1)\quad d_1=4,\;\;d_2d_3>167,\qquad 2)\quad d_1=5,
\;\;d_2d_3>73\qquad 3)\quad  d_1\ge 6.\label{b26}
\eea
A brief analysis of (\ref{b26}) and restriction $d_3<d_1d_2-d_1-d_2$ give a short list of non-symmetric semigroups $S_3$ with $d_1=4,5$, where 
inequality (\ref{b22}) might be broken. Below we present Table 1 of 13 semigroups $S_3$, where $\rho_3>1.1547$.
%%%%%%%%%%%%%%%%%%%%%%%%%%%%%%%%
\begin{center}
Table 1
\end{center}
%%%%%%%%%%%%%%%\
$$
\begin{array}{|c||c|c|c|c|}\hline
 S_3 &\rho_3 & e_3 & w_3 & h_0  \\\hline\hline
\quad\la 4,5,7\ra\quad &\quad 1.35225\quad &\quad 1.07143\quad 
&\quad 0.57143\quad &\quad 0.52273\quad \\\hline
\la 4,7,9\ra  &1.25988 &0.86363 & 0.54545 &0.51667 \\\hline
\la 4,5,11\ra &1.34840 &1.18750 & 0.625 & 0.53704 \\\hline
\la 4,7,13\ra &1.25794 &1.04545 & 0.63636 & 0.54412 
\\\hline
\la 4,7,17\ra &1.28338 &0.96428 & 0.57143 & 0.52439  \\\hline
\la 4,9,11\ra &1.20605 &0.76667 & 0.53333 &0.51316 \\
\hline
\la 4,9,15\ra &1.20493 &0.9     & 0.6   & 0.53571 \\\hline
\la 4,11,13\ra &1.17074&0.71053 & 0.52632 & 0.51087 
\\\hline\hline
\la 5,6,7\ra &1.24212  &0.85    & 0.6    & 0.53704 \\\hline
\la 5,6,8\ra &1.22644  &0.9     & 0.6    & 0.53571 \\\hline
\la 5,7,8\ra &1.19523  &0.79167 & 0.58333 &0.53226 \\\hline
\la 5,7,9\ra &1.18322  &0.71429 & 0.57143 &0.52941 \\\hline
\la 5,8,9\ra &1.15950  &0.80769 & 0.61539 & 0.54412\\\hline\hline
\end{array}
$$
%%%%%%%%%%%%%%%%%%%%%%%%%%%%%%%%%%%%%%%%
\vspace{.5cm}

Ten semigroups $S_3$ among 13 in Table 1,
$$
\la 4,7,9\ra,\;\la 4,7,17\ra,\;\la 4,9,11\ra,\;\la 4,9,15\ra,\;
\la 4,11,13\ra,\;\la 5,6,7\ra,\;\la 5,6,8\ra,\;\la 5,7,8\ra,\;
\la 5,7,9\ra,\;\la 5,8,9\ra,
$$
have $\rho_3>1.1547$ and $e_3<1$, and, according to inequality (\ref{b20ā}), provide $w_3<2/3$. The rest three semigroups $S_3$,
\bea
\la 4,5,7\ra,\quad\la 4,5,11\ra,\quad\la 4,7,13\ra,\nonumber
\eea
have $\rho_3>1.1547$ and $e_3>1$, but a direct calculation gives
\bea
w_3(\la 4,5,7\ra)=0.5714,\qquad
w_3(\la 4,5,11\ra)=0.625,\qquad
w_3(\la 4,7,13\ra)=0.6363.\nonumber
\eea
Thus, Theorem is proven.$\;\;\;\;\;\;\;\;\Box$
\end{proof}
%%%%%%%%%%%%%%%%%%%%%%%%%%%%%%%%%%%%%%%%%%%%%%%%%
\subsection{Lower and upper bounds of $h_r=K_r/g_3^{r+1}$, $r\ge 1$.}\label{l24}
%%%%%%%%%%%%%%%%%%%%%%%%%%%%%%%%%%%%%%%%%%%%%%%%%
Consider four syzygy identities (\ref{b2}),
\bea
Y_1-X_1=0,\qquad Y_2-X_2=2\pi_3,\qquad Y_3-X_3=
6\pi_3K_0,\qquad Y_4-X_4=24\pi_3K_1.\label{b27}
\eea
Substituting into the fourth of them polynomial relations $Y_4(Y_1,Y_2)$ and $X_4(X_1,X_2,X_3$ (see \cite{fe21}, formulas (37)),
\bea
2Y_4=Y_2^2+2Y_1^2Y_2-Y_1^4,\qquad
6X_4=X_1^4-6X_1^2X_2+8X_1X_3+3X_2^2,\nonumber
\eea
we obtain
\bea
3(Y_2^2+2Y_1^2Y_2-Y_1^4)=X_1^4-6X_1^2X_2
+8X_1X_3+3X_2^2+144\pi_3K_1.\label{b28}
\eea
Continuing to substitute into (\ref{b28}) three first identities (\ref{b27}) and one more $2Y_3=3Y_1Y_2-Y_1^3$,
\bea
Y_2=Y_1^2-4K_0Y_1+12K_1+\pi_3,\nonumber
\eea
and combining it with identity (\ref{b6}), we arrive at equality in rescaled variables $h_1,h_0,u,v$,
\bea
12h_1=2u(2h_0-1)+2-v,\qquad h_1=\frac{K_1}{g_3^2}.\label{b29}
\eea
Making use of bounds (\ref{b8}, \ref{b14}, \ref{b19}) for $u,v,h_0$, we arrive at lower and upper bounds for $h_1$,
\bea
0.13889\simeq\frac{5}{36}<h_1<\frac{11}{54}\simeq 0.20370.\label{b30}
\eea

Next, supplement (\ref{b27}) with one more identity from (\ref{b2}),
\bea
Y_5-X_5=60\pi_3K_2,\label{b31}
\eea
and substitute there polynomial relations $Y_5(Y_1,Y_2)$ and $X_5(X_1,X_2,X_3)$ (see \cite{fe21}, formulas (37)), and the other three first identities (\ref{b2}). Skipping lengthy calculations, we present the final equality (see \cite{fe21}, formulas (39)),
\bea
Y_1^3-2K_0Y_1^2+4\pi_3K_0+24K_2=Y_2(Y_1+2K_0),\nonumber
\eea
and combine it with identity (\ref{b6}). We present the final formula in rescaled variables $h_2,h_0,u,
v$,
\bea
12h_2=(2h_0-1)u(u-1)+2h_0(1-v),\qquad h_2=\frac{K_2}{g_3^3}.\label{b32}
\eea
Making use of bounds (\ref{b8}, \ref{b14}, \ref{b19}) for $u,v,h_0$, we arrive at lower and upper bounds for $h_2$,
\bea
0.05555\simeq\frac{1}{18}<h_2<\frac{1}{9}\simeq 0.11111.\label{b33}
\eea

The ratios $h_r\!=\!K_r/g_3^{r+1}$, $r\!\ge\! 0$, will be referred to as rescaled genera of numerical semigroup $S_3$.

To study the bounds of $h_r$, $r\ge 3$, we make worth of Theorem 2 in \cite{fe21}, applied to semigroups $S_3$: there exists an algebraic equation $R(h_0,h_1,h_2,h_r)\!=\!0$, where the polynomial  
$R(t_1,t_2,t_3,t_4)$ is irreducible over a ring $A[t_1,t_2,t_3,t_4]$. Avoiding lengthy formulas of $R(h_0,h_1,h_2,h_r)$ with growing $r$, we present one algebraic equation for $r=3$, making use of formula (42) in \cite{fe21},
\bea
&&\left(10h_3-18h_1^2+vh_0^2-\frac{v^2}{24}\right)\Delta_1=\Delta_2^2,\label{b34}\\
&&\Delta_1=3h_1-2h_0^2+\frac{v}{4},\quad\Delta_2=6h_2-6h_0h_1+\frac{v}{2}h_0.\nonumber
\eea
Substituting (\ref{b29},\ref{b32}) into (\ref{b34}), we obtain
\bea
\Delta_1=\frac1{2}(2h_0-1)(u-1-2h_0),\qquad\Delta_2=u\Delta_1,\nonumber
\eea
that substantially simplifies equality (\ref{b34}),
\bea
10h_3=18h_1^2-vh_0^2+\frac{v^2}{24}+u^2\Delta_1.\label{b35}
\eea
Making use of bounds (\ref{b8},\ref{b14},\ref{b19},\ref{b30}), we arrive at lower and upper bounds for $h_3$,
\bea
0.02443<h_3<0.07515.\label{b36}
\eea
%%%%%%%%%%%%%%%%%%%%%%%%$%%%%%%%%%%%%%%%%%%%%%%%%
\subsection{Rescaled genera $h_r$ in symmetric semigroups $S_3$}\label{l25}
%%%%%%%%%%%%%%%%%%%%%%%%%%%%%%%%%%%%
Every symmetric semigroup $S_3$ is a complete intersection \cite{he70} and therefore we can apply formulas (67) in \cite{fe21} to calculate the six first rescaled genera $h_r$, $1\le r\le 6$, and find their lower and upper bounds. Bearing in mind that
\bea
h_0=\frac1{2},\quad u=1,\quad 0<v\le\frac1{4},\nonumber
\eea
we get
\bea
&&\hspace{-2cm}
h_1=\;\frac1{6}\left(1-\frac1{2}v\right),
\hspace{4cm} 
0.14583\simeq\frac{7}{48}<h_1<
\;\frac1{6}\simeq 0.16666,\label{b37}\\
&&\hspace{-2cm}
h_2=\frac1{12}\;(1-v),\hspace{4.5cm}0.06250\simeq\frac1{16}<h_2<\frac1{12}
\simeq 0.08333,\nonumber\\
&&\hspace{-2cm}
h_3=\frac1{20}\left(1-\frac{3}{2}v+\frac1{3}v^2\right),
\hspace{2.7cm}
0.03229\simeq\frac{31}{960}< h_3<\frac1{20}\simeq 0.05,\nonumber\\
&&\hspace{-2cm}
h_4=\frac1{30}\;\left(1-2v+v^2\right),
\hspace{3.2cm}
0.01875\simeq\frac{3}{160}<h_4<\frac1{30}\simeq 0.03333,\nonumber\\
&&\hspace{-2cm}
h_5=\frac1{42}\left(1-\frac{5}{2}v+2v^2-\frac1{4}v^3\right),
\hspace{1.4cm}
0.01181\simeq\frac{127}{10752}< h_5<\frac1{42}
\simeq 0.02381,\nonumber\\
&&\hspace{-2cm}
h_6=\frac1{56}\left(1-3v+\frac{10}{3}v^2-v^3\right),\hspace{1.6cm}
0.00791\simeq\frac{85}{10752}<h_6<\frac1{56}\simeq 0.01785.\nonumber
\eea

Making use of (\ref{b19}) and comparing the upper and lower bounds in (\ref{b30},\ref{b33},\ref{b36}) with those in (\ref{b37}), we conclude that the domain of variation of $h_r$, $0\le r\le 3$, in non-symmetric 
semigroups $S_3$ contains their upper and lower bounds in symmetric semigroups $S_3$,
\bea
0.5=0.5=\!\!&h_0&\!\!=0.5<0.55555,\label{b38}\\
0.13889<0.14593<\!\!&h_1&\!\!<0.16666<0.20370,\nonumber\\
0.05555<0.06250<\!\!&h_2&\!\!<0.08333<0.11111,\nonumber\\
0.02443<0.03229<\!\!&h_3&\!\!<0.05000<0.07515.\nonumber
\eea
e.g., $\;h_1(\la 4,7,13\ra)=0.169406$, $\;h_1(\la 5,8,9\ra)=0.169694$.

We leave open a question of whether such a property for the other rescaled genera $h_r$ with arbitrary $r\ge 4$ exists in numerical semigroups $S_3$. 
%%%%%%%%%%%%%%%%%%%%%%%%%%%%%%%%%%%%%%%%%%%%%%%%%%%%%%%%%%%

\end{document}